\documentclass{amsart}

\newcommand{\supp}{\mbox{supp}\,}
\renewcommand{\span}{\mbox{span}\,}
\renewcommand{\r}{{\Bbb R}}
\newcommand{\z} {{\Bbb Z}}
\newcommand{\n} {{\Bbb N}}
\newcommand{\ddd}{,\dots,}
\newcommand{\lll}{\left(}
\newcommand{\rrr}{\right)}
\newcommand{\ex}[1]{e^{2\pi i{#1}}}

\newcommand{\cD}{{\mathcal D}}

\newcommand{\bN}{{\mathbb N}}

\newcommand{\bQ}{{\mathbb Q}}
\newcommand{\bR}{{\mathbb R}}
\newcommand{\bC}{{\mathbb C}}

\newcommand{\be}{\begin{equation}}
\newcommand{\ee}{\end{equation}}
\newcommand{\ba}{\begin{eqnarray}}
\newcommand{\ea}{\end{eqnarray}}
\newcommand{\ban}{\begin{eqnarray*}}
\newcommand{\ean}{\end{eqnarray*}}

\newtheorem{theorem}{Theorem}[section]
\newtheorem{lemma}[theorem]{Lemma}

\theoremstyle{definition}
\newtheorem{definition}[theorem]{Definition}
\newtheorem{example}[theorem]{Example}

\theoremstyle{remark}
\newtheorem{remark}[theorem]{Remark}

\theoremstyle{proposition}
\newtheorem{proposition}[theorem]{Proposition}

\theoremstyle{proposition}
\newtheorem{corollary}[theorem]{Corollary}

\numberwithin{equation}{section}

\begin{document}

\title[$p$-Adic multiresolution analysis and wavelet frames]
{$p$-Adic multiresolution analysis\\and wavelet frames}

\author{S.~Albeverio}
\address{Universit\"at Bonn, Institut f\"ur Angewandte
Mathematik, Abteilung Sto\-chas\-tik, Wegelerstra\ss e 6, D-53115 Bonn
and Interdisziplinäres Zentrum f\"ur Komplexe Systeme, Universit\"at Bonn,
R\"omerstra\ss e 164 D-53117, Bonn, Germany}
\email{albeverio@uni-bonn.de}
\thanks{The first and the third authors were supported in part by
DFG Project 436 RUS 113/809. The second author was supported in part by Grants
06-01-00471 and 07-01-00485 of RFBR.
The third author was supported in part by Grant
06-01-00457 of RFBR}

\author{S.~Evdokimov}
\address{St.-Petersburg Department of Steklov Institute  of Mathematics, St.-Petersburg,
 Fontanka-27, 191023 St. Petersburg, RUSSIA }
 \email{evdokim@pdmi.ras.ru}

\author{M.~Skopina}
\address{Department of Applied Mathematics and Control Processes,
St. Petersburg State University, \ Universitetskii pr.-35,
198504 St. Petersburg, Russia.}
\email{skopina@MS1167.spb.edu}


\date{}



\begin{abstract}
We study $p$-adic multiresolution analyses (MRAs). A complete
characterisation of  test   functions generating  MRAs (scaling
functions) is given. We prove that only  $1$-periodic test functions
may be taken as orthogonal scaling functions. We also suggest a method for
the construction of wavelet functions and prove that any wavelet
function generates a  $p$-adic wavelet frame.
\end{abstract}

\maketitle

\section{Introduction}
\label{s1}

In the early nineties a general scheme for the
construction of wavelets (of real argument) was developed. This
scheme is based on the notion of multiresolution analysis (MRA in the sequel)
introduced by Y.~Meyer and S.~Mallat~\cite{Mallat-1}, \cite{Meyer-1} (see also,
e.g.,  ~\cite{31}, ~\cite{NPS}). Immediately specialists started to
implement new wavelet systems.
Nowadays it is difficult to find an engineering area where wavelets are
not applied.

 In the $p$-adic setting, the situation is as follows.
In 2002 S.~V.~Kozyrev~\cite{Koz0} found a compactly supported
$p$-adic wavelet basis for ${ L}^2(\bQ_p)$ which is an analog of the Haar basis.
It turned out that these wavelets were eigenfunctions of
$p$-adic pseudo-differential operators~\cite{Koz2}.
J.J.~Benedetto and R.L.~Benedetto~\cite{Ben-Ben} conjectured that other $p$-adic wavelets
with the same set of translations can not be constructed because this set is not a group,
and  the corresponding MRA-theory can not be developed.
Another conjecture was rised by  A.~Khrennikov and V.~Shelkovich~\cite{Kh-Sh1}.
They assumed that the equality
\begin{equation}
\label{62.0-3}
\phi(x)=\sum_{r=0}^{p-1}\phi\Big(\frac{1}{p}x-\frac{r}{p}\Big),
\quad x\in \bQ_p,
\end{equation}
may be considered as a {\it refinement equation} for the Haar MRA
generating Kozyrev's wavelets. A solution $\phi$ to this equation
({\it a refinable function}) is the characteristic function of the
unit disc. We note that equation (\ref{62.0-3}) reflects a {\it
natural} ``self-similarity'' of the space $\bQ_p$: the unit disc
$B_{0}(0)=\{x: |x|_p \le 1\}$ is represented as the union
$\bigcup_{r=0}^{p-1}B_{-1}(r)$ of $p$ mutually {\it disjoint} discs
$B_{-1}(r)=\bigl\{x: |x-r|_p \le p^{-1}\bigr\}$
(see~\cite[I.3, Examples 1,2.]{Vl-V-Z}).
Following this idea, the notion of $p$-adic MRA
was introduced and a general scheme for its
construction was described in~\cite{S-Sk-1}. Also, using
(\ref{62.0-3}) as a generating refinement equation, this scheme was
realized to construct the $2$-adic Haar MRA. In contrast to the real
setting, the {\it refinable function} $\phi$ generating the Haar MRA
is {\em periodic}, which implies the existence of {\em infinitly many
different} orthonormal wavelet bases in the same Haar MRA. One of
them coincides with Kozyrev's wavelet basis.
The authors of~\cite{Kh-Sh-S} described a wide class of  functions generating
a MRA, but all of these functions are $1$-periodic.
In the present paper we prove that there exist no other orthogonal
test scaling functions generating
a MRA, except for those described in~\cite{S-Sk-1}.
Also, the MRAs generated by arbitrary test scaling functions (not necessary
orthogonal) are considered. We thoroughly study these scaling functions
and develop a method
to construct a wavelet frame based on a given MRA.

Here and in what follows, we shall systematically use the
notation and the results from~\cite{Vl-V-Z}.
Let $\bN$, $\z$, $\r$, $\bC$ be the sets of positive integers, integers,
real numbers, complex numbers, respectively.
The field $\bQ_p$ of $p$-adic numbers is defined as the completion
of the field of rational numbers $\bQ$ with respect to the
non-Archimedean $p$-adic norm $|\cdot|_p$. This $p$-adic norm
is defined as follows: $|0|_p=0$; if $x\ne 0$, $x=p^{\gamma}\frac{m}{n}$,
where $\gamma=\gamma(x)\in \z$
and the integers $m$, $n$ are not divisible by $p$, then
$|x|_p=p^{-\gamma}$.
The norm $|\cdot|_p$ satisfies the strong triangle inequality
$|x+y|_p\le \max(|x|_p,|y|_p)$.
The canonical form of any $p$-adic number $x\ne 0$ is
\begin{equation}
\label{2}
x=p^{\gamma}(x_0 + x_1p + x_2p^2 + \cdots),
\end{equation}
where $\gamma=\gamma(x)\in \z$, \ $x_j\in D_p:=\{0,1,\dots,p-1\}$, $x_0\ne 0$,
$j=0,1,\dots$. We shall write the $p$-adic numbers
$k=k_{0}+k_{1}p+\cdots+k_{s-1}p^{s-1}$, $k_j\in D_p$, $j=0,1,\dots,s-1$,
following the usual  form, as in  the real analysis: $k=0,1,\dots,p^s-1$.

Denote by $B_{\gamma}(a)=\{x\in \bQ_p: |x-a|_p \le p^{\gamma}\}$
the disc of radius $p^{\gamma}$ with the center at a point $a\in \bQ_p$,
$\gamma \in \z$.
Any two balls in $\bQ_p$ either are disjoint or one
contains the other.

There exists the Haar measure $dx$ on $\bQ_p$ which is  positive,
invariant under the shifts, i.e., $d(x+a)=dx$, and normalized by
$\int_{|\xi|_p\le 1}\,dx=1$.
A complex-valued function $f$ defined on $\bQ_p$ is called
{\it locally-constant} if for any $x\in \bQ_p$ there exists
an integer $l(x)\in \z$ such that
$f(x+y)=f(x)$, $y\in B_{l(x)}(0)$.
Denote by ${\cD}$ the linear space of locally-constant compactly
supported functions (so-called test functions)~\cite[VI.1.,2.]{Vl-V-Z}.
The space ${\cD}$ is an analog of the  Schwartz space in the real analysis.

The Fourier transform of $\varphi\in {\cD}$ is defined as
$$
{\widehat\phi}(\xi)=F[\varphi](\xi)=\int_{\bQ_p}\chi_p(\xi\cdot x)\varphi(x)\,dx,
\ \ \ \xi \in \bQ_p,
$$
where $\chi_p(\xi\cdot x)=e^{2\pi i\{\xi x\}_p}$ is the additive character for
the field $\bQ_p$, $\{\cdot\}_p$ is a fractional part of a number $x\in \bQ_p$.
The Fourier transform is a linear isomorphism taking ${\cD}$ into
${\cD}$. The Fourier transform is extended to ${ L}^2(\bQ_p)$ in a
standard way. If $f\in{ L}^2(\bQ_p)$, $0\ne a\in \bQ_p$, \ $b\in \bQ_p$,
then~\cite[VII,(3.3)]{Vl-V-Z}:
\begin{equation}
\label{014}
F[f(ax+b)](\xi)
=|a|_p^{-1}\chi_p\Big(-\frac{b}{a}\xi\Big)F[f(x)]\Big(\frac{\xi}{a}\Big).
\end{equation}
According to~\cite[IV,(3.1)]{Vl-V-Z},
\begin{equation}
\label{14.1}
F[\Omega(p^{-k}|\cdot|_p)](x)=p^{k}\Omega(p^k|x|_p), \quad k\in \z,
\quad x \in \bQ_p,
\end{equation}
where $\Omega(t)=1$ for $t\in [0,\,1]$; $\Omega(t)=0$ for $t\not\in [0,\,1]$.

\section{Multiresolution analysis}
\label{s2}

Let us consider the set
$$
I_p=\{a=p^{-\gamma}\big(a_{0}+a_{1}p+\cdots+a_{\gamma-1}p^{\gamma-1}\big):
\gamma\in \bN; a_j\in D_p; j=0,1,\dots,\gamma-1\}.
$$
It is well known that
$\bQ_p=B_{0}(0)\cup\cup_{\gamma=1}^{\infty}S_{\gamma}$, where
$S_{\gamma}=\{x\in \bQ_p: |x|_p = p^{\gamma}\}$. Due to
(\ref{2}), $x\in S_{\gamma}$, $\gamma\ge 1$, if and only if
$x=x_{-\gamma}p^{-\gamma}+x_{-\gamma+1}p^{-\gamma+1}+\cdots+x_{-1}p^{-1}+\xi$,
where $x_{-\gamma}\ne 0$, $\xi \in B_{0}(0)$. Since
$x_{-\gamma}p^{-\gamma}+x_{-\gamma+1}p^{-\gamma+1}
+\cdots+x_{-1}p^{-1}\in I_p$, we have a ``natural'' decomposition of
$\bQ_p$ into a union of mutually  disjoint discs:
$\bQ_p=\bigcup_{a\in I_p}B_{0}(a)$.
So, $I_p$ is a {\em ``natural'' set of shifts} for $\bQ_p$.

\begin{definition}
\label{de1} \rm
A collection of closed spaces
$V_j\subset L^2(\bQ_p)$, $j\in\z$, is called a
{\it multiresolution analysis {\rm(}MRA{\rm)} in $ L^2(\bQ_p)$} if the
following axioms hold

(a) $V_j\subset V_{j+1}$ for all $j\in\z$;

(b) $\bigcup_{j\in\z}V_j$ is dense in $ L^2(\bQ_p)$;

(c) $\bigcap_{j\in\z}V_j=\{0\}$;

(d) $f(\cdot)\in V_j \Longleftrightarrow f(p^{-1}\cdot)\in V_{j+1}$
for all $j\in\z$;

(e) there exists a function $\phi \in V_0$
such that $V_0:=\overline{\span\{\phi(\cdot-a),\ a\in I_p\}}$.
\end{definition}

The function $\phi$ from axiom (e) is called {\em scaling}.
One also says that a MRA is generated by its scaling function $\phi$
(or $\phi$ generates the MRA).
It follows  immediately from axioms (d) and (e) that
\be
V_j:=\overline{\span\{\phi(p^{-j}x-a),\ a\in I_p\}},\quad j\in \z.
\label{17}
\ee

An important class of MRAs consists of those
generated by so-called {\em orthogonal scaling functions}.
A scaling function  $\phi$ is said to be orthogonal if
 $\{\phi(\cdot-a), a\in I_p\}$ is an orthonormal
basis for $V_0$. Consider such a MRA.
Evidently, the functions $p^{j/2}\phi(p^{-j}\cdot-a)$,
$a\in I_p$, form an orthonormal basis for $V_j$, $j\in\z$.
According to the standard scheme (see, e.g.,~\cite[\S 1.3]{NPS})
for the construction of MRA-based wavelets, for each $j$, we define
a space $W_j$ ({\em wavelet space}) as the orthogonal complement
of $V_j$ in $V_{j+1}$, i.e., $V_{j+1}=V_j\oplus W_j$, $j\in \z$,
where $W_j\perp V_j$, $j\in \z$. It is not difficult to see that
\begin{equation}
\label{61.0}
f\in W_j \Longleftrightarrow f(p^{-1}\cdot)\in W_{j+1},
\quad\text{for all}\quad j\in \z
\end{equation}
and $W_j\perp W_k$, $j\ne k$.
Taking into account axioms (b) and (c), we obtain
\begin{equation}
\label{61.1}
{\bigoplus\limits_{j\in\z}W_j}= L^2(\bQ_p)
\quad \text{(orthogonal direct sum)}.
\end{equation}
If we now find   functions $\psi^{(\nu)} \in W_0$, $\nu\in A$,
such that the functions $\psi^{(\nu)}(x-a)$, $a\in~I_p, \nu\in A$, form an orthonormal
basis for $W_0$, then, due to~(\ref{61.0}) and (\ref{61.1}),
the system $\{p^{j/2}\psi^{(\nu)}(p^{-j}\cdot-a), a\in I_p, j\in\z , \nu\in A\}$
is an orthonormal basis for $ L^2(\bQ_p)$.
Such a function $\psi$ is called a {\em wavelet function} and
the basis is a {\em wavelet basis}.

Another interesting class of scaling functions consists of
functions $\phi$  so that  $\{\phi(\cdot-a), a\in I_p\}$ is a Riesz system.
Probably, adopting the ideas developed for the real setting, one can
use MRAs generated by such functions  $\phi$ for construction of
dual biorthogonal wavelet systems. This topic is, however, out of our consideration
in the present paper.

In  Section~\ref{s3} we will discuss how to construct a $p$-adic wavelet frame
based on an arbitrary MRA generated by a test function.

Let $\phi$ be an orthogonal scaling function for a MRA $\{V_j\}_{j\in\z}$. Since
the system $\{p^{1/2}\phi(p^{-1}\cdot-a), a\in I_p\}$
is a basis for $V_1$ in this case, it follows from axiom (a) that
\begin{equation}
\label{62.0-2*}
\phi=\sum_{a\in I_p}\alpha_a\phi(p^{-1}\cdot-a),
\quad \alpha_a\in \bC.
\end{equation}
We see that the function $\phi$ is a solution of a
special kind of functional equation. Such equations are called
{\em refinement equations}, and their solutions are called  {\em refinable functions}
\footnote{Usually the terms ``refinable function'' and ``scaling function'' are
synonyms in the literature, and they are used in both senses: as a solution
to the refinable equation and as a function generating MRA.
We separate here the meanings of these terms.}.
It will be shown in Section~\ref{s3} that any test scaling function
(not necessary orthogonal) is  refinable.

A natural way for the construction of a MRA (see, e.g.,~\cite[\S 1.2]{NPS})
is the following. We start with a refinable function $\phi$
and define the spaces $V_j$ by~(\ref{17}).
It is clear that axioms (d) and (e) of Definition~\ref{de1} are fulfilled.
Of course, not any such function $\phi$ provides axiom $(a)$.
In the  real setting, the relation $V_0\subset V_{1}$ holds
if and only if the refinable function satisfies a refinement equation.
The situation is different in the $p$-adic case.. Generally speaking, a refinement
equation (\ref{62.0-2*})  does not imply the including property
$V_0\subset V_{1}$ because the set of shifts $I_p$  does not
form a group.
Indeed, we need all the functions $\phi(\cdot-b)$,
$b\in I_p$, to belong to the space $V_1$, i.e., the identities
$\phi(x-b)=\sum_{a\in I_p}\alpha_{a,b}\phi(p^{-1}x-a)$ should be
fulfilled for all $b\in I_p$. Since $p^{-1}b+a$ is not in $I_p$ in general,
we  can not state that $\phi(x-b)$ belongs to  $V_1$
for all $b\in I_p$.
Nevertheless, we will see below that a wide class of refinable equations
provide the including property.

Providing axiom (a) is a key moment for the construction of MRA.
Axioms (b) and (c) are fulfilled for a wide class of functions $\phi$
because of the following statements.

\begin{theorem}
\label{th1-2*}
If $\phi \in  L^2(\bQ_p)$ and $\widehat\phi$ is compactly supported,
then axiom $(c)$ of Definition~{\rm\ref{de1}} holds for the spaces
$V_j$ defined by~(\ref{17}).
\end{theorem}

\begin{proof} Let $\widehat\phi\subset B_M(0)$, $M\in\z$.
Assume that a function $f\in L^2(\bQ_p)$ belongs to any space $V_j$, $j\in\z$.
 Given  $j\in\n$ and $\epsilon>0$, there exists a  function
 $f_\epsilon:=\sum_{a\in I_p}\alpha_a\phi(p^j\cdot-a)$, where the sum is finite,
  such that $\|f-f_\epsilon\|<\epsilon$. Using~(\ref{014}),
  it is not difficult to see that
  $\supp\widehat f_\epsilon\subset\supp\widehat\phi(p^{-j}\cdot)$,
  which yields that   $\widehat f_\epsilon(\xi)=0$ for any
  $\xi\not\in B_{M-j}(0)$.  Due to the Plancherel theorem,
  it follows that $\widehat f=0$ almost everywhere
  on $B_{M-j}(0)$. Since $j$ is an arbitrary positive integer,
 $\widehat f$ is equivalent to zero on $Q_p$.
\end{proof}

Another sufficient condition for axiom (c) was given in~\cite{Kh-Sh-S}:

\begin{theorem}
\label{th1-4*}
If $\phi \in  L^2(\bQ_p)$ and the system $\{\phi(x-a):a\in I_p\}$ is orthonormal,
then axiom $(c)$ of Definition~{\rm\ref{de1}} holds for the spaces
$V_j$ defined by~(\ref{17}).
\end{theorem}

\begin{theorem}
Let $\phi \in  L^2(\bQ_p)$, the spaces $V_j$,
$j\in \z$, be defined by~(\ref{17}),
and let $\phi(\cdot-b)\in\cup_{j\in \z} V_j$ for any $b\in Q_p$.
Axiom $(b)$ of Definition~{\rm\ref{de1}} holds for the spaces
$V_j$, $j\in\z$,   if and only if
\be
\bigcup\limits_{j\in\z}{\rm supp\,}\widehat\phi(p^{j}\cdot)=\bQ_p.
\label{dnn14}
\ee
\label{th1-3*}
\end{theorem}

\begin{remark}
It is not difficult to see that the assumption
$\phi(\cdot-b)\in\cup_{j\in \z} V_j$ for any $b\in Q_p$
is fulfilled whenever $\phi$ is a refinable function
and $\widehat\phi\subset B_0(0)$. We will see that
this assumption is also valid for a wide class of refinable
functions  $\phi$ for which $\widehat\phi\not\subset B_0(0)$.
\end{remark}

\begin{proof}
First of all we show that the space
$\overline {\cup_{j\in {\z}} V_j}$ is invariant with respect to all
shifts. Let $f\in \cup_{j\in {\z}} V_j$, $b\in\bQ_p$. Evidently,
$\phi(p^{-k}{\cdot}-t)\in\cup_{j\in \z} V_j$ for any $t\in Q_p$
and for any $k\in\z$.
Since the $L_2$-norm is invariant with
respect to the shifts,  it follows  that
$f(\cdot -b) \in \overline{\cup_{j\in\z} V_j}$.
 If now $g\in\overline{\cup_{j\in\z} V_j}$, then
approximating $g$ by the functions  $f\in \cup _{j\in\z} V_j$,
again using  the invariance of $L_2$-norm with respect to the shifts
, we derive $g(\cdot -b) \in \overline{\cup_{j\in\z} V_j}$.

For $X\subset L^2(\bQ_p)$, set $\widehat X=\{\ \widehat f: f\in X\}$.
By the Wiener theorem for $L_2$ (see, e.g., \cite{NPS}; all the arguments
of the proof given there may be repeated word for word with replacing
 $\bR$ by $\bQ_p$), a closed subspace
$X$ of the space $L^2(\bQ_p)$ is invariant with respect to the shifts
if and only if $\widehat X=L_2(\Omega)$ for some set
$\Omega\subset\bQ_p$. If now $X=\overline{\cup_{j\in\z}V_j}$, then
$\widehat X=L_2(\Omega)$. Thus $X=L^2(\bQ_p)$ if and only if $\Omega=\bQ_p$. Set
$\phi_j=\phi(p^{-j}\cdot),\ \ \Omega_0=\cup_{j\in\z}{\rm supp}\, \widehat\phi_j$
and prove that $\Omega=\Omega_0$. Since $\phi_j\in V_j,$ $j\in\z$,
we have ${\rm supp}\, \widehat\phi_j\subset\Omega$, and hence
$\Omega_0 \subset \Omega$.
Now assume that $\Omega\backslash\Omega_0$ contains a set of
positive measure $\Omega_1$. Let $f\in V_j$. Given   $\epsilon>0$,
there exists a  function
 $f_\epsilon:=\sum_{a\in I_p}\alpha_a\phi(p^j\cdot-a)$, where the sum is finite,
  such that $\|f-f_\epsilon\|<\epsilon$. Using~(\ref{014}),
  we see that $\supp\widehat f_\epsilon\subset\supp\widehat\phi(p^{-j}\cdot)$,
  which yields that   $\widehat f_\epsilon(\xi)=0$ for any
  $\xi\not\in \Omega_1$.  Due to the Plancherel theorem,
  it follows that $\widehat f=0$ almost everywhere  on $\Omega_1.$
Hence the same is true for any $f\in \cup _{j\in\z} V_j$.
Passing to the limit we deduce that that the Fourier transform
of any $f\in X$ is equal to zero almost everywhere on  $\Omega_1$, i.e.,
$L_2(\Omega)=L_2(\Omega_0)$.
It remains to note that
${\rm supp\,} \widehat\phi_j={\rm supp\,}\widehat\phi({p}^{j}\cdot) $
\end{proof}

A real analog of Theorem~\ref{th1-3*} was proved by
 C.~de~Boor,  R.~DeVore and A.~Ron in~\cite5.

\section{Refinable functions}
\label{s3}

We are going to study $p$-adic refinable functions $\phi$. Let us restrict
ourselves to the consideration of  $\phi\in {\cD}$.
Evidently, each $\phi\in {\cD}$ is a
$p^M$-periodic  function for some $M\in \z$.
Denote by ${\cD}_N^M$ the set of all $p^M$-periodic  functions supported on $B_N(0)$.
Taking the Fourier transform
of the equality $\phi(x-p^M)=\phi(x)$, we obtain
$\chi_p(p^M\xi)\widehat\phi(\xi)=\widehat\phi(\xi)$, which holds for all $\xi$
if and only if $\supp\widehat\phi\subset B_M(0)$. Thus, the set ${\cD}_N^M$
consists of all locally constant functions $\phi$ such that $\supp\phi\subset B_N(0)$,
 $\supp\widehat\phi\subset B_M(0)$.

\begin{proposition}
Let $\phi,\psi \in  L^2(\bQ_p)$, $\supp \phi,\supp \psi\subset B_N(0)$, $N\ge0$,
 and let $b\in I_p$, $|b|_p\le p^N$. If
\be
\psi(\cdot-b)\in \overline{\span\{\phi(p^{-1}\cdot-a),\ a\in I_p\}}
\label{19}
\ee
then
\be
\psi(x-b)=\sum_{k=0}^{p^{N+1}-1}h^\psi_{k,b}\phi\Big(\frac{x}{p}-\frac{k}{p^{N+1}}\Big)
\ \ \ \forall x\in Q_p.
\label{18}
\ee
\label{p1}
\end{proposition}

\begin{proof}
Given  $\epsilon>0$, there exist  functions
 $$
 f_\epsilon:=\sum_{a\in I_p\atop |a|_p\le p^{N+1}}\alpha_a\phi(p^j\cdot-a),\ \ \
g_\epsilon:=\sum_{a\in I_p\atop |a|_p> p^{N+1}}\alpha_a\phi(p^j\cdot-a),
$$
  where the sums are finite,
  such that $\|\psi(\cdot-b)-f_\epsilon-g_\epsilon\|<\epsilon$.
  If $x\in B_N(0)$, $|a|_p>p^{N+1}$, then
$|p^{-1}x-a|_p>p^{N+1}$ and hence $\phi(p^{-1}x-a)=0$. So, $g_\epsilon(x)=0$
whenever $x\in B_N(0)$. If $x\not\in B_N(0)$, then $\phi(x-b)=0$ and
$\phi(p^{-1}x-a)=0$ for all $a\in I_p$, $|a|_p\le p^{N+1}$.
So, $\phi(\cdot-b)-f_\epsilon(x)=0$ whenever $x\not\in B_N(0)$.
It follows that
$$
\|\psi(\cdot-b)-f_\epsilon\|^2=\int\limits_{B_N(0)}|\psi(\cdot-b)-f_\epsilon|^2=
\int\limits_{B_N(0)}|\psi(\cdot-b)-f_\epsilon-g_\epsilon|^2\le \epsilon^2.
$$
Hence
$$
\psi(\cdot-b)\in
\overline{\span\{\phi(p^{-1}\cdot-a),\ a\in I_p,\ |a|_p\le p^{N+1}\}},
$$
which implies~(\ref{18}).
\end{proof}

\begin{corollary}
If $\phi \in  L^2(\bQ_p)$ is a refinable  function and
 $\supp \phi\subset B_N(0)$, $N\ge0$,  then  its refinement equation is
\begin{equation}
\label{62.0-5}
\phi(x)=\sum_{k=0}^{p^{N+1}-1}h_{k}\phi\Big(\frac{x}{p}-\frac{k}{p^{N+1}}\Big)
\ \ \ \forall x\in Q_p.
\end{equation}
\label{c3}
\end{corollary}

The proof immediately follows from Proposition~\ref{p1}.

\begin{corollary}
Let $\phi \in  L^2(\bQ_p)$ be a scaling  function of a MRA.
If  $\supp \phi\subset B_N(0)$, $N\ge0$,  then  $\phi$ is a refinable
function satisfying~(\ref{62.0-5}).
\label{c1}
\end{corollary}

The proof follows by combining axiom $(a)$ of Definition~{\rm\ref{de1}} with
Proposition~\ref{p1}.

Taking the Fourier transform of~(\ref{62.0-5})  and using
(\ref{014}),  we can rewrite the refinable equation in the form
\begin{equation}
\label{62.0-6}
{\widehat\phi}(\xi)=m_0\Big(\frac{\xi}{p^{N}}\Big){\widehat\phi}(p\xi),
\end{equation}
where
\begin{equation}
\label{62.0-7-1}
m_0(\xi)=\frac{1}{p}\sum_{k=0}^{p^{N+1}-1}h_{k}\chi_p(k\xi)
\end{equation}
is a trigonometric polynomial. It is clear that $m_0(0)=1$ whenever $\widehat\phi(0)\ne0$.

\begin{proposition}
\label{pr1-1}
If $\phi\in  L^2(\bQ_p)$ is a solution of refinable
equation~(\ref{62.0-5}), ${\widehat\phi}(0)\ne 0$, ${\widehat\phi}(\xi)$ is continuous at
the point $0$, then
\begin{equation}
\label{62.0-8}
{\widehat\phi}(\xi)={\widehat\phi}(0)\prod_{j=0}^{\infty}m_0\Big(\frac{\xi}{p^{N-j}}\Big).
\end{equation}
\end{proposition}

\begin{proof}
Since~(\ref{62.0-5}) implies~(\ref{62.0-6}), after
iterating {\rm(\ref{62.0-6})} $J$ times, $J\ge 1$, we have
$$
{\widehat\phi}(\xi)=\prod_{j=0}^{J}m_0\Big(\frac{\xi}{p^{N-j}}\Big)
{\widehat\phi}(p^{J}\xi).
$$
Taking into account that ${\widehat\phi}(\xi)$ is continuous at
the point $0$ and the fact that $|p^{N}\xi|_p=p^{-N}|\xi|_p\to 0$
as $N\to +\infty$ for any $\xi\in\bQ_p$, we obtain {\rm(\ref{62.0-8})}.
\end{proof}

\begin{corollary}
If $\phi \in  {\cD}_N^M$ is a refinable function, $N\ge0$, and ${\widehat\phi}(0)\ne 0$, then
(\ref{62.0-8}) holds.
\label{c2}
\end{corollary}

This statement follows immediately from Corollary~\ref{c1} and
Proposition~\ref{pr1-1}.

\begin{lemma}
Let
$
{\widehat\phi}(\xi)=C\prod_{j=0}^{\infty}
m_0\Big(\frac{\xi}{p^{N-j}}\Big),
$
where $m_0$ is a trigonometric polynomial with $m_0(0)=1$
and $C\in\r$.
If $\supp\widehat\phi\subset B_M(0)$,
then there exist at least
$\lll p^{M+N}-\frac{\deg m_0}{p-1}\rrr$ integers $n$ such that
$0\le n<p^{M+N}$ and $\widehat\phi\lll\frac{n}{p^M}\rrr=0$.
\label{l1}
\end{lemma}

\begin{proof}
First of all we note that $\widehat\phi$ is a $p^N$-periodic
function satisfying~(\ref{62.0-6}).
Denote by $O_p$ the set of positive integers not divisible by~$p$.
Since $\supp\widehat\phi\subset B_M(0)$,
we have $\widehat\phi\lll\frac{k}{p^{M+1}}\rrr=0$ for all
$k\in O_p$
. By the definition of $\widehat\phi$ the equality
$\widehat\phi\lll\frac{k}{p^{M+1}}\rrr=0$ holds if and only if
there exists $\nu=1-N\ddd M+1$ such that $m_0\Big(\frac{k}{p^{N+\nu}}\Big)=0$.
Set
$$
\sigma_\nu:=\left\{l\in O_p:\ l<p^{N+\nu},
m_0\Big(\frac{l}{p^{N+\nu}}\Big)=0,\
m_0\Big(\frac{l}{p^{N+\mu}}\Big)\ne0\ \forall \mu=1-N\ddd\nu-1\right\},
$$
$v_\nu:=\sharp\,\sigma_\nu$. Evidently, $\sigma_\nu\subset O_p^{\prime}$
for all $\nu$, where $O_p^{\prime}=\{k\in O_p:\ k<p^{M+N+1}\}$,
and $\sigma_{\nu^{\prime}}\cap\sigma_{\nu}=\emptyset$ whenever
$\nu^{\prime}\ne\nu$. If $\widehat\phi\lll\frac{k}{p^{M+1}}\rrr=0$ for some
$k\in O_p$, then there exist a unique $\nu=1-N\ddd M+1$
and a unique $l\in\sigma_\nu$ such that $k\equiv l\pmod{p^{N+\nu}}$. Moreover,
for any $l\in\sigma_\nu$ there are exactly $p^{M-\nu+1}$ integers $k\in O_p^{\prime}$
(including~$l$) satisfying the above comparison.
It follows that
\be
\sum\limits_{\nu=1-N}^{M+1} p^{M-\nu+1}v_\nu=\sharp\,O_p^{\prime}=p^{M+N}(p-1).
\label{10}
\ee
Now if $l\in\sigma_\nu$, $\nu\le M$, then
$\widehat\phi\lll\frac{p^\gamma k}{p^{M}}\rrr=0$
 for all $\gamma=0,1\ddd M-\nu$, $k=l+rp^{N+\nu}$, $r=0,1\ddd p^{M-\nu-\gamma}-1$,
 i.e., each $l\in\sigma_\nu$ generates at least $1+p+\dots+p^{M-\nu}$
 distinct positive integers $n<p^{M+N}$ for which
 $\widehat\phi\lll\frac{n}{p^{M}}\rrr=0$. Hence
 \ban
v:= \sharp\,\left\{n:\ n=0,1\ddd p^{M+N}-1,
\widehat\phi\lll\frac{n}{p^{M}}\rrr=0\right\}\ge
\\
\sum\limits_{\nu=1-N}^{M}(1+p+\dots+p^{M-\nu})v_\nu=
\frac{1}{p-1}\sum\limits_{\nu=1-N}^{M}(p^{M-\nu+1}-1)v_\nu=
\\
\frac{1}{p-1}\sum\limits_{\nu=1-N}^{M+1}(p^{M-\nu+1}-1)v_\nu.
\ean
Since $\sum\limits_{\nu=1-N}^{M+1}v_\nu\le \deg m_0$, by using~(\ref{10}),
we obtain
$$
v\ge\frac{1}{p-1}\lll\sum\limits_{\nu=1-N}^{M+1}p^{M-\nu+1}v_\nu-
\deg m_0\rrr\ge p^{M+N}-\frac{\deg m_0}{p-1}.
$$
\end{proof}

\begin{theorem}
\label{t1}
Let $\phi\in{\cD}_N^M$, $N\ge0$ and ${\widehat\phi}(0)\ne 0$. If
\be
\phi(\cdot-b)\in \overline{\span\{\phi(p^{-1}\cdot-a),\ a\in I_p\}}
\label{27}
\ee
for all $b\in I_p$, $|b|_p\le p^N$, then there exist at least
$p^{M+N}-p^N$ integers $l$ such that $0\le l<p^{M+N}$ and
$\widehat\phi\lll\frac{l}{p^{M}}\rrr=0$.
\end{theorem}

\begin{proof}
Let $b\in I_p$, $|b|_p\le p^N$. Because of Proposition~\ref{p1},
we can rewrite~(\ref{27}) in the form
$$
\phi(x-b)=\sum_{k=0}^{p^{N+1}-1}h_{k,b}\phi\Big(\frac{x}{p}-\frac{k}{p^{N+1}}\Big)
\ \ \ \forall x\in \bQ_p.
$$
Taking the Fourier transform, we obtain
\begin{equation}
\label{12}
{\widehat\phi}(\xi)\chi_p(b\xi)=m_b\Big(\frac{\xi}{p^{N}}\Big){\widehat\phi}(p\xi),
\ \ \ \forall \xi\in \bQ_p,
\end{equation}
where $m_b$ is a trigonometric polynomial, $\deg m_b<p^{N+1}$.
Combining~(\ref{12}) for $b=0$ with~(\ref{12}) for arbitrary $b$,
we obtain
$$
{\widehat\phi}(p\xi)\lll m_0\Big(\frac{\xi}{p^{N}}\Big)\chi_p(b\xi)-
 m_b\Big(\frac{\xi}{p^{N}}\Big)\rrr=0
\ \ \ \forall \xi\in \bQ_p,
$$
which is equivalent to
\be
F(\xi):={\widehat\phi}(p^{N+1}\xi)\lll m_0(\xi)\chi_p(p^{N}b\xi)-
 m_b(\xi)\rrr=0
\ \ \ \forall \xi\in \bQ_p.
\label{13}
\ee
Since  $\supp F\subset B_{M+N+1}(0)$ and $F$ is a $1$-periodic function,
(\ref{13}) holds if and only if $\widehat\phi\lll\frac{l}{p^{M+N+1}}\rrr=0$,
$l=0,1\ddd p^{M+N+1}-1$.

First suppose that $\deg m_0\ge p^N(p-1)$, i.e.,
$$
m_0(\xi)=\sum\limits_{k=0}^{K}h_k\chi_p(k\xi),\ \ \  h_K\ne0,
$$
where $K=K_Np^N+K_{N-1}p^{N-1}+\dots+K_0$,
$K_j\in D_p$, $j=0,1\ddd N$, $K_N=p-1$ (indeed,
if $K_N<p-1$, then $\deg m_0= K\le (p-2)p^N+(p-1)(1+p+\dots+p^{N-1})=
p^{N+1}-p^N-1<p^N(p-1)$). Set $b:= p-p^{-N}K$. It is not difficult to see that
$b\in I_p$, $|b|_p\le p^N$ and $K+bp^N=p^{N+1}$. We see that
the degree of the polynomial $t(\xi):=m_0(\xi)\chi_p(p^{N}b\xi)- m_b(\xi)$
is exactly $p^{N+1}$, and hence there exist at most $p^{N+1}$ integers $l$
such that $0\le l<p^{M+N+1}$, $t\lll\frac{l}{p^{M+N+1}}\rrr=0$. Thus,
$$
\sharp\,\left\{l:\ l=0,1\ddd p^{M+N+1}-1,
\widehat\phi\lll\frac{l}{p^{M}}\rrr=0\right\}\ge p^{M+N+1}-p^{N+1}.
$$
Taking into account that $\widehat\phi$ is a $p^N$-periodic function,
we obtain
\be
 \sharp\,\left\{l:\ l=0,1\ddd p^{M+N}-1,
\widehat\phi\lll\frac{l}{p^{M}}\rrr=0\right\}\ge p^{M+N}-p^{N}.
\label{14}
\ee
It remains to note that~(\ref{14}) is also fulfilled  whenever
$\deg m_0<p^N(p-1)$ because of Lemma~\ref{l1} and Corollary~\ref{c2}.
\end{proof}

\begin{theorem}
\label{t2}
Let $\phi,\psi\in{\cD}_N^M$, $N\ge0$, ${\widehat\phi}(0)\ne 0$, and let
\be
\psi(\cdot)\in \overline{\span\{\phi(p^{-1}\cdot-a),\ a\in I_p\}}
\label{28}
\ee
If there exist at least
$p^{M+N}-p^N$ integers $l$ such that $0\le l<p^{M+N}$ and
$\widehat\phi\lll\frac{l}{p^{M}}\rrr=0$, then
\be
\psi(x-b)=\sum_{a\in I_p}\alpha^\psi_{a,b}\phi(p^{-1}x-a) \ \ \
\forall b\in \bQ_p,
\label{11}
\ee
where the sum  is finite. In particular. if $\phi$ is a refinable function,
then
\be
\phi(x-b)=\sum_{a\in I_p}\alpha_{a,b}\phi(p^{-1}x-a) \ \ \
\forall b\in \bQ_p.
\label{11}
\ee

\end{theorem}

\begin{proof}
First we assume that $b\in Q_p$, $|b|_p\le p^{N}$, $b\ne0$, and prove that
\be
\psi(x-b)=\sum_{k=0}^{p^{N+1}-1}g_{k,b}\phi\Big(\frac{x}{p}-\frac{k}{p^{N+1}}\Big)
\ \ \ \forall x\in \bQ_p.
\label{29}
\ee
Because of Proposition~\ref{p1}, we have~(\ref{29}) for $b=0$.
Taking the Fourier transform of~(\ref{29}), we obtain
\begin{equation}
\label{12}
{\widehat\psi}(\xi)\chi_p(b\xi)=n_b\Big(\frac{\xi}{p^{N}}\Big){\widehat\phi}(p\xi),
\ \ \ \forall \xi\in \bQ_p,
\end{equation}
where $n_b$ is a trigonometric polynomial, $\deg n_b<p^{N+1}$.
Substituting~(\ref{12}) for $b=0$,
we reduce~(\ref{12}) for arbitrary $b$ to
$$
{\widehat\phi}(p\xi)\lll n_0\Big(\frac{\xi}{p^{N}}\Big)\chi_p(b\xi)-
 n_b\Big(\frac{\xi}{p^{N}}\Big)\rrr=0
\ \ \ \forall \xi\in \bQ_p,
$$
which is equivalent to
\be
F(\xi):={\widehat\phi}(p^{N+1}\xi)\lll n_0(\xi)\chi_p(p^{N}b\xi)-
 n_b(\xi)\rrr=0
\ \ \ \forall \xi\in \bQ_p.
\label{13}
\ee
Since  $\supp F\subset B_{M+N+1}(0)$ and $F$ is a $1$-periodic function,
(\ref{13}) is equivalent to
$$
F\lll\frac{l}{p^{M+N+1}}\rrr=0,
\forall l=0,1\ddd p^{M+N+1}-1,
$$
which holds if and only if
\ba
n_b\lll\frac{l}{p^{M+N+1}}\rrr=
n_0\lll\frac{l}{p^{M+N+1}}\rrr\chi_p\lll\frac{bl}{p^{M+1}}\rrr,\hspace{2cm}
\label{15}
\ea
for all $l=0,1\ddd p^{M+N+1}-1$ such that
$\widehat\phi\lll\frac{l}{p^{M}}\rrr\ne0$. Because of  $p^M$-periodicity of
$\widehat\phi$, there exist at least $p(p^{M+N}-p^N)$ integers
$l=0,1\ddd p^{M+N+1}-1$ such that
$\widehat\phi\lll\frac{l}{p^{M}}\rrr=0$. So, we can find $n_b$ by  solving
the linear system~(\ref{15}) with respect to the unknown coefficients
of $n_b$. Taking the Fourier transform of~(\ref{12}), we obtain
\be
\psi(x-b)=\sum_{k=0}^{p^{N+1}-1}g_{k,b}\phi\Big(\frac{x}{p}-\frac{k}{p^{N+1}}\Big)
\ \ \ \forall x\in \bQ_p.
\label{16}
\ee
Next let  $b\in Q_p$, $|b|_p= p^{N+1}$, i.e., $b=b_{N+1}p^{N+1}+b^\prime$,
$b_{N+1}\in D_p$, $b_{N+1}\ne0$, $|b^\prime|_p\le p^{N}$.
Using~(\ref{16}) with $b=b^\prime$, we have
$$
\psi(x-b)=\sum_{k=0}^{p^{N+1}-1}g_{k,b^\prime}
\phi\Big(\frac{x}{p}-\frac{k}{p^{N+1}}-\frac{b_{N+1}}{p^{N+2}}\Big)=
\sum_{k=0}^{p^{N+1}-1}g_{k,b^\prime}
\phi\Big(\frac{x}{p}-\frac{pk+b_{N+1}}{p^{N+2}}\Big).
$$
Taking into account that
$$
pk+b_{N+1}\le p(p^{N+1}-1)+(p-1)=p^{N+2}-1,
$$
we derive
$$
\psi(x-b)=\sum_{k=0}^{p^{N+2}-1}g_{k,b}\phi\Big(\frac{x}{p}-\frac{k}{p^{N+2}}\Big)
\ \ \ \forall x\in \bQ_p.
$$
Similarly, we can prove by induction on $n$ that
$$
\psi(x-b)=\sum_{k=0}^{p^{N+n+1}-1}g_{k,b}\phi\Big(\frac{x}{p}-\frac{k}{p^{N+n+1}}\Big)
\ \ \ \forall x\in \bQ_p,
$$
whenever $b\in \bQ_p$, $|b|_p= p^{N+n}$.
\end{proof}

\begin{theorem}
\label{t3}
A function $\phi\in{\cD}_N^M$, $N\ge0$,  with
${\widehat\phi}(0)\ne 0$ generates a MRA if and only if

(1) $\phi$ is refinable;

(2) there exist at least
$p^{M+N}-p^N$ integers $l$ such that $0\le l<p^{M+N}$ and
$\widehat\phi\lll\frac{l}{p^{M}}\rrr=0$.
\end{theorem}

\begin{proof}
If $\phi$ is a scaling function of a MRA, then (1) follows from
Corollary~\ref{c1} and (2) follows from (1) and Theorem~\ref{t1}.

Now let conditions (1), (2) be fulfilled. Define the spaces $V_j$, $j\in\z$,
by~(\ref{17}). Axioms (d) and (e), evidently, hold.
Axiom (a) follows from Theorem~\ref{t2}. Axiom (b) follows from
Theorems~\ref{t2} and ~\ref{th1-3*}. Axiom (c) follows from
Theorems~\ref{th1-2*}.
\end{proof}

\begin{example}
Let $p=2$, $N=2$, $M=1$  $\phi$ be defined by (\ref{62.0-8}),
where $\widehat\phi(0)\ne0$, $m_0$ is given by (\ref{62.0-7-1}),
$m_0(1/4)=m_0(3/8)=m_0(7/16)=m_0(15/16)=0$ and $m_0(0)=1$.
It is not difficult to see that ${\rm supp\,}\widehat\phi\subset B_1(0)$,
${\rm supp\,}\widehat\phi\not\subset B_0(0)$ and
$\widehat\phi\Big(\frac{1}{2}\Big)=\widehat\phi\Big(\frac{3}{2}\Big)=
\widehat\phi\Big(\frac{5}{2}\Big)=
\widehat\phi(1)=0$, i.e, all the assumptions of Theorem~\ref{t3}
are fulfilled.
\label{e1}
\end{example}

\begin{remark}
The above example is typical. Similarly, taking into account
the arguments of the proof of Lemma~\ref{l1}, one can easily construct a lot of
functions $\phi$ generating a MRA for arbitrary  $p$, $M>0$ and large enough
 $N$. Moreover, it is possible to provide $\deg m_0\le2^N$.
 \label{r1}
\end{remark}

\section{Orthogonal scaling functions}
\label{s4}

Now we are going to describe all orthogonal
scaling functions  $\phi\in{\cD}_N^M$.

\begin{theorem}
\label{th1-5*}
Let  $\phi\in{\cD}_N^M$,  $M,N\ge0$.
If    $\{\phi(x-a):a\in I_p\}$
is an orthonormal system, then
\be
\sum_{l=0}^{p^{M+N}-1}\left|{\widehat\phi}\lll\frac{l}{p^M}\rrr\right|^2
\chi_p\lll\frac{lk}{p^{M+N}}\rrr=p^N\delta_{k 0},
\quad k=0,{1},\dots,{p^{N}-1}.
\label{20}
\ee
\end{theorem}

\begin{proof}
Let $a\in I_p$. Due to the orthonormality of $\{\phi(x-a):a\in I_p\}$,
using the Plancherel theorem, we have
$$
\delta_{a 0}=
\langle\phi(\cdot),\phi(\cdot-a)\rangle
\int\limits_{\bQ_p}\phi(x)\overline{\phi(x-a)}\,dx
=\int\limits_{B_M(0)}|{\widehat\phi}(\xi)|^2\chi_p(a\xi)\,d\xi.
$$
Let $\xi\in B_M(0)$. There exists a unique $l=0,1\ddd p^{M+N}-1$
such that  $\xi\in B_{-N}\lll b_l\rrr$, $b_l=\frac{l}{p^M}$. It follows that
\ban
\int\limits_{B_M(0)}|{\widehat\phi}(\xi)|^2\chi_p(a\xi)\,d\xi
=\sum_{k=0}^{p^{M+N}-1}\int\limits_{|\xi-b_l|_p\le p^{-N}}
|{\widehat\phi}(\xi)|^2\chi_p(a\xi)\,d\xi\qquad\qquad\qquad\qquad
\\
=\sum\limits_{l=0}^{p^{M+N}-1}|{\widehat\phi}(b_l)|^2
\int\limits_{|\xi-b_l|_p\le p^{-N}}\chi_p(a\xi)\,d\xi
=\sum\limits_{l=0}^{p^{M+N}-1}|{\widehat\phi}(b_l)|^2\chi_p(ab_l)
\int\limits_{|\xi|_p\le p^{-N}}\chi_p(a\xi)\,d\xi
\\
=\frac{1}{p^{N}}\Omega(|p^{N}a|_p)
\sum_{l=0}^{p^{M+N}-1}|{\widehat\phi}(b_l)|^2\chi_p(ab_l).
\ean

To prove~(\ref{20}) it only remains to note that
 $\Omega(|p^{N}a|_p)= 0$ whenever $a\in I_p$,
 $p^Na\ne0,1\ddd p^N-1$.
\end{proof}

\begin{lemma}
Let $c_0,\ddd c_{n-1}$ be mutually distinct elements of the unit
circle $\{z\in\bC:\ |z|=1\}$. Suppose that there exist nonzero
reals $x_j$, $j=0,1\ddd n-1$, such that
\be
\sum_{j=0}^{n-1}c_j^kx_j=\delta_{k 0},\ \ \ k=0,1\ddd n-1.
\label{21}
\ee
Then $x_j=1/n$ for all $j$, and up to reordering
\be
c_j=c_0\ex{j/n},\ \ \ j=0,1\ddd n-1.
\label{22}
\ee
\label{l2}
\end{lemma}

\begin{proof}
In accordance with  Cramer's rule we have $x_j=\frac{\Delta_j}{\Delta}$,
$0\le j\le n-1$, where $\Delta=V(c)$ is the Vandermonde determinant
corresponding to $c=(c_0\ddd c_{N-1})$, and $\Delta_j$ is obtained from
$\Delta$ by replacing the $j$-th column with the transpose of the row
$(1,0\ddd 0)$. A straightforward computation shows that
$$
\Delta_j=(-1)^jV(c^{(j)})\prod\limits_{k\ne j}c_k,
$$
where $c^{(j)}$ is obtained from $c$ by removing the $j$-th coordinate. Thus,
\ba
x_j=(-1)^j\frac{V(c^{(j)})}{V(c)}\prod\limits_{k\ne j}c_k
=(-1)^j\prod\limits_{k\ne j}c_k\prod\limits_{k>l\atop k,l\ne j}
({c_k-c_l})\Big/{\prod\limits_{k>l}(c_k-c_l)}
\nonumber
\\
\prod\limits_{k\ne j}\frac{c_k}{c_k-c_j}=
\prod\limits_{k\ne j}\frac{1}{1-c^{-1}_kc_j}.
\label{35}
\ea
Next, for any $\alpha\in\r$, we have
\ban
1-e^{i\alpha}=2\sin\frac{\alpha}{2}
\lll\sin\frac{\alpha}{2}-i\cos\frac{\alpha}{2}\rrr=
2\sin\frac{\alpha}{2}e^{i\lll\frac{\alpha}{2}-\frac{\pi}{2}\rrr}.
\ean
Let us define $\alpha_j$, $j=0,1\ddd n-1$, by $c_j=e^{i\alpha_j}$.
Then from the above arguments and~(\ref{35}) it follows that
$$
x_j=\prod\limits_{k\ne j}\frac{1}{1-c^{-1}_kc_j}=
 e^{i\gamma}\sum\limits_{k\ne j}\lll2\sin\frac{\alpha_k-\alpha_j}{2}\rrr^{-1},
$$
where
$$
\gamma=\sum\limits_{k\ne j}\frac{\alpha_k-\alpha_j+\pi}{2}=
\theta-\frac{n}{2}\alpha_j,\ \ \ \theta=
\frac12\lll(n-1)\pi+\sum\limits_{k=0}^{n-1}\alpha_k\rrr
$$
By the lemma's hypothesis $x_j\in\r$, whence
$\gamma\equiv0\pmod{\pi}$ and consequently
$n\alpha_j\equiv2\theta\pmod{2\pi}$. Thus up to
reordering $\alpha_j=\alpha_0+\frac{2\pi j}{n}$,
which implies~(\ref{22}), and consequently that $x_j=1/n$ for all $j$.
\end{proof}

\begin{theorem}
\label{t4}
Let  $\phi\in{\cD}_N^M$ be an orthogonal
scaling function and  $\widehat\phi(0)\ne0$.
Then ${\rm supp\,{\widehat\phi}}\subset B_0(0)$.
\end{theorem}

\begin{proof}
Without loss of generality, we can assume that  $M,N\ge0$.
Combining Theorems~\ref{t3} and \ref{th1-5*}, we have
$$
\sum_{j=0}^{p^{N}-1}\left|{\widehat\phi}\lll\frac{l_j}{p^M}\rrr\right|^2
\chi_p\lll\frac{l_jk}{p^{M+N}}\rrr=p^N\delta_{k 0},
\quad k=0,{1},\dots,{p^{N}-1}.
$$
By Lemma~\ref{l2}, $l_j=l_0+jp^M$ and ${\widehat\phi}\lll\frac{l_j}{p^M}\rrr=1$.
Taking into account that  $\widehat\phi(0)\ne0$, we deduce
$l_0=0$, i.e., ${\widehat\phi}(j)=1$, $j=0,1\ddd p^{N}-1$.
Since $\widehat\phi$ is a $p^N$-periodic function,
it follows from Theorem~\ref{t3}  that
$\widehat\phi\lll\frac{l}{p^M}\rrr=0$ for all $l\in\z$
not divisible by $p^M$.
This yields ${\rm supp\,{\widehat\phi}}\subset B_0(0)$.
\end{proof}

So any test function $\phi$ generating a MRA belongs to the class
${\cD}_N^0$. All such functions were described in~\cite{Kh-Sh-S}.
The following theorem summarizes these results.

\begin{theorem}
\label{th1-6*}
Let ${\widehat\phi}$ be defined by {\rm(\ref{62.0-8})}, where $m_0$ is
the trigonometric polynomial~{\rm(\ref{62.0-7-1})} with $m_0(0)=1$.
If $m_0\big(\frac{k}{p^{N+1}}\big)=0$ for all $k=1,\dots,p^{N+1}-1$
not divisible by $p$, then  $\phi\in{\cD}_N^0$.
If, furthermore,
$\big|m_0\big(\frac{k}{p^{{N+1}}}\big)\big|=1$ for all $k=1,\dots,p^{N+1}-1$
divisible by $p$, then $\{\phi(x-a):a\in I_p\}$ is an orthonormal system.
Conversely, if  ${\rm supp\,{\widehat\phi}}\subset B_0(0)$
and the system $\{\phi(x-a):a\in I_p\}$
is orthonormal, then $\big|m_0\big(\frac{k}{p^{{N+1}}}\big)\big|=0$
whenever $k$ is not divisible
by $p$, and $\big|m_0\big(\frac{k}{p^{{N+1}}}\big)\big|=1$ whenever $k$
is divisible by $p$, \ $k=1,2,\dots,p^{N+1}-1$.
\end{theorem}

\section{Construction of wavelet frames}
\label{s4}

\begin{definition}
Let $H$ be a Hilbert space. A system $\{f_n\}_{n=1}^\infty\subset H$
is said to be a frame if there exist positive constants $A,B$
({\em frame boundaries}) such that
$$
A\|f\|^2\le \sum_{n=1}^\infty|\langle f,f_n\rangle|^2\le B\|f\|^2
\ \ \ \forall f\in H.
$$
\end{definition}

We are interested in the construction of $p$-adic wavelet frames, i.e.,
frames in $L_2(\bQ_p)$ consisting of functions $p^{j/2}\psi^{(\nu)}(p^{-j}\cdot-a)$,
$a\in I_p$, $\nu\in A$, where $A$ is a  finite set.

We will restrict ourselves to the consideration of the case $p=2$.

Our general scheme of construction looks as follows. Let $\{V_j\}_{j\in\z}$
be a MRA. As above, we define
the wavelet space $W_j$, $j\in \z$,  as the orthogonal complement
of $V_j$ in $V_{j+1}$, i.e., $V_{j+1}=V_j\oplus W_j$.
 It is not difficult to see that
$f\in W_j$ if and only if $f(2^{j}\cdot)\in W_0$,
and $W_j\perp W_k$ whenever $j\ne k$. If now there exists a
function $\psi\in L_2(Q_2)$ ({\em wavelet function}) such that
\be
W_0=\overline{\span\{\psi(\cdot-a),\ a\in I_2\}},
\label{23}
\ee
then we have a wavelet system
$\{2^{j/2}\psi(2^{-j}\cdot-a), a\in I_2, j\in\z \}$.
It will be proved that such a system is a frame   in $L_2(Q_2)$ whenever
 $\psi$ is compactly supported.

\begin{theorem}
Let $\{V_j\}_{j\in\z}$ be a MRA, $\psi$ be a wavelet function.
If $\psi$ is compactly supported, then the corresponding
wavelet system
$\{2^{j/2}\psi(2^{-j}\cdot-a), a\in I_2, j\in\z \}$
 is a frame in $L_2(Q_2)$.
\label{t5}
 \end{theorem}

\begin{proof}
First we will prove that  the system $\{\psi(\cdot-a),\ a\in I_2\}$
is a frame in  the wavelet space $W_0$.
Let  $\supp\psi\subset B_N(0)$, $N\ge0$. Set
\ban
W_0^0&=&{\span\{\psi(\cdot-a),\ a\in I_2, |a|_2\le2^N\}},
\\
W_0^n&=&{\span\{\psi(\cdot-a),\ a\in I_2, |a|_2=2^{N+n}\}},\ \ \ n\in\n.
\ean
It is not difficult to see that the spaces $W_0^n$, $n=0,1,\dots$,
are mutually orthogonal. Each function $f\in W_0$ may be represented
in the form $f=f^0+f^1+\cdots$, where $f^0=f\Big|_{B_N(0)}$,
$f^n=f\Big|_{B_{N+n}(0)\setminus B_{N+n-1}(0)}$, $n\in \n$.
Due to~(\ref{23}), given $\epsilon>0$, there exists a finite
sum $\sum\limits_{a\in I_2}\alpha_a\psi(\cdot-a)=:f_\epsilon$ such that
$\|f-f_\epsilon\|<\epsilon$.  If $|x|_2\le2^N$, then
$f_\epsilon(x)=
\sum\limits_{a\in I_2\atop |a|_2\le2^N}\alpha_a\psi(x-a)=:f^0_\epsilon(x)$.
Since $\supp f^0\subset B_N(0)$, $\supp f_\epsilon^0\subset B_N(0)$, we have
$$
\|f-f_\epsilon\|^2\ge\int\limits_{B_N(0)}|f-f_\epsilon|^2
=\int\limits_{B_N(0)}|f^0-f^0_\epsilon|^2=\|f^0-f^0_\epsilon\|^2.
$$
It follows that $f^0\in W_0^0$. Similarly, $f^n\in W_0^n$, $n\in\n$.
Thus we proved that
\be
W_0=W_0^0\oplus W_0^1\oplus W_0^2\oplus\dots.
\label{26}
\ee
Since $W_0^0$ is a finite dimensional space and
$\{\psi(\cdot-a),\ a\in I_2, |a|_2\le2^N\}$ is a representing system for
$W_0^0$, this system is a frame. Hence there exist
 positive constants $A,B$ such that
$$
 A\|f^0\|^2\le\sum_{a\in I_2}|\langle f^0,\psi(\cdot-a)\rangle|^2\le B\|f^0\|^2
  \ \ \ \forall f\in W_0^0.
$$
 If  $f^1\in W^1_0$, we have
\ban
\sum_{a\in I_2\atop |a|_2=2^{N+1}}|\langle f^1,\psi(\cdot-a)\rangle|^2=
\sum_{a\in I_2\atop |a|_2\le2^N}|\langle f^1,\psi(\cdot-a-2^{-N-1})\rangle|^2=
\\
\sum_{a\in I_2\atop |a|_2\le2^N}|\langle f^1(\cdot+2^{-N-1}),\psi(\cdot-a)\rangle|^2\ge
 A\|f^1(\cdot+2^{-N-1})\|^2=A\|f^1\|^2.
\ean
Let now $f^n\in W_0^n$, $n>1$.  Then
\ban
\sum_{a\in I_2\atop |a|_2=2^{N+n}}
\left|\left\langle f^n,\psi(\cdot-a)\right\rangle\right|^2=
\sum\limits_{k=0}^{2^{n-1}-1}
\sum_{a\in I_2\atop |a|_2\le2^N}
\left|\left\langle f^n,\psi\lll\cdot-a-\frac{2k+1}{2^{N+n}})
\rrr\right\rangle\right|^2=
\\
\sum\limits_{k=0}^{2^{n-1}-1}
\sum_{a\in I_2\atop |a|_2\le2^N}
\left|\left\langle f^n\lll\cdot+\frac{2k+1}{2^{N+n}}\rrr
\Omega(|2^N\cdot|),\psi\lll\cdot-a)
\rrr\right\rangle\right|^2\ge
\\
 A\sum\limits_{k=0}^{2^{n-1}-1}
 \left\|f^n\lll\cdot+\frac{2k+1}{2^{N+n}}\rrr\Omega(|2^N\cdot|)\right\|^2=
  A\sum\limits_{k=0}^{2^{n-1}-1}
\left\|f^n\Omega\lll\left|2^N\lll\cdot-\frac{2k+1}{2^{N+n}}\rrr\right|\rrr\right\|^2=
\\
  A\sum\limits_{k=0}^{2^{n-1}-1}
\left\|f^n\Big|_{B_N\lll\frac{2k+1}{2^{N+n}}\rrr}\right\|^2=
A\|f^n\|^2.
\ean
Taking into account~(\ref{26}), we  derive
\be
 A\|f\|^2\le\sum_{a\in I_2}|\langle f,\psi(\cdot-a)\rangle|^2
 \ \ \ \forall f\in W_0.
\label{25}
\ee
 Similarly we can prove the upper frame estimation
\be
 \sum_{a\in I_2}|\langle f,\psi(\cdot-a)\rangle|^2\le B\|f\|^2
 \ \ \ \forall f\in W_0.
\label{24}
\ee
Combining~(\ref{25}) with~(\ref{24}), we
deduce that the system $\{\psi(\cdot-a),\ a\in I_2\}$
is a frame in   $W_0$. Evidently,
the system $\{2^{j/2}\psi(2^{-j}\cdot-a), a\in I_2, \}$
is a frame in   $W_j$ with the same frame boundaries for any $j\in\z$.
Since ${\bigoplus\limits_{j\in\z}W_j}= L^2(\bQ_2)$, it follows that
the union of these frames is a frame in $L_2(\bQ_2)$.
\end{proof}

Now we discuss how to construct a desirable wavelet function $\psi$.
Let a MRA $\{V_j\}_{j\in\z}$  is generated by a scaling function $\phi\in{\cD}_N^M$.
First of all we should provide $\psi\in V_1$. Let us look for $\psi$
in the form
$$
\psi(x)=\sum_{k=0}^{2^{N+1}-1}g_{k}\phi\Big(\frac{x}{2}-\frac{k}{2^{N+1}}\Big)
$$
Taking the Fourier transform of~(\ref{62.0-5})  and using
(\ref{014}),  we have
$$
{\widehat\psi}(\xi)=n_0\Big(\frac{\xi}{2^{N}}\Big){\widehat\phi}(2\xi),\ \ \
$$
where $n_0$ is a trigonometric polynomial ({\em wavelet mask}) given by
$$
n_0(\xi)=\frac{1}{2}\sum_{k=0}^{2^{N+1}-1}g_{k}\chi_2(k\xi)
$$
Evidently, $\psi\in{\cD}_N^M$.
By Theorem~\ref{t1}, there exist at least
$2^{M+N}-2^N$ integers $l$ such that $0\le l<2^{M+N}$,
$\widehat\phi\lll\frac{l}{2^{M}}\rrr=0$. Choose $n_0$
satisfying the following property: if
$\widehat\phi\lll\frac{l}{2^{M}}\rrr\ne0$ for some
$l=0,1\ddd 2^{M+N}-1$, then $n_0\lll\frac{l}{2^{M+N}}\rrr=0$.
This yields that $\widehat\psi\lll\frac{l}{2^{M}}\rrr=0$ whenever
$\widehat\phi\lll\frac{l}{2^{M}}\rrr\ne0$, $0\le l<2^{M+N}$.

Let $a, b\in I_2$.
Using the Plancherel theorem and the arguments of Theorem~\ref{t1}, we have
\ban
\langle\phi(\cdot-a),\psi(\cdot-b)\rangle=
\int\limits_{\bQ_2}\phi(x-a)\overline{\psi(x-b)}\,dx
=\int\limits_{B_M(0)}{\widehat\phi}(\xi)
\overline{\widehat\psi(\xi)}\chi_2((b-a)\xi)\,d\xi=
\\
\sum_{k=0}^{2^{M+N}-1}\int\limits_{|\xi-2^{-M}l|_2\le 2^{-N}}
{\widehat\phi}(\xi)
\overline{\widehat\psi(\xi)}\chi_2((b-a)\xi)\,d\xi=\qquad\qquad\qquad\qquad
\\
\sum\limits_{l=0}^{2^{M+N}-1}{\widehat\phi}\lll\frac{l}{2^M}\rrr
\overline{\widehat\psi\lll\frac{l}{2^M}\rrr}
\int\limits_{|\xi-2^{-M}l|_2\le 2^{-N}}\chi_2(a\xi)\,d\xi=0.
\ean
It follows that $\overline{\span\{\psi(\cdot-a),\ a\in I_2\}}\perp V_0$.
On the other hand, due to Theorem~\ref{t2}, we have
$\overline{\span\{\psi(\cdot-a),\ a\in I_2\}}\subset V_1$.
Hence,
\be
\overline{\span\{\psi(\cdot-a),\ a\in I_2\}}\subset W_0.
\label{30}
\ee
It is clear from that proof of Theorem~\ref{t2} that
\ba
\psi\lll x-\frac{l}{p^{N}}\rrr&=&
\sum_{k=0}^{p^{N+1}-1}g_{kl}\phi\Big(\frac{x}{p}-\frac{k}{p^{N+1}}\Big),
\ \ \ \ l=0,1\ddd 2^N-1,
\label{31}
\\
\phi\lll x-\frac{l}{p^{N}}\rrr&=&
\sum_{k=0}^{p^{N+1}-1}h_{kl}\phi\Big(\frac{x}{p}-\frac{k}{p^{N+1}}\Big),
\ \ \ l=0,1\ddd 2^N-1.
\label{32}
\ea
Consider these equalities as a linear system with respect to the
unknowns $X_k:=\phi\Big(\frac{x}{p}-\frac{k}{p^{N+1}}\Big)$, $k=0,1\ddd 2^{N+1}-1$.
If the system~(\ref{31}),~(\ref{32})  has a solution, then
$$
\span\left\{\phi\Big(\frac{\cdot}{p}-a\Big),
\ a\in I_2, |a|_2\le2^{N+1}\right\}
\subset{\span\{\psi(\cdot-a),\ a\in I_2, |a|_2\le2^{N}\}}.
$$
This evidently implies
$W_0\subset\overline{\span\{\psi(\cdot-a),\ a\in I_2\}}.$
Taking into account~(\ref{30}), we deduce that $\psi$ is a wavelet function.

It is not quite clear whether  the system~(\ref{31}),~(\ref{32})  has a solution for
arbitrary $\phi$ and $\psi$, but we will show how to succeed in the case
$\deg m_0\le2^N$. The construction of such masks can easily be done (see
Example~\ref{e1} and Remark~\ref{r1}).

Assume that $\deg m_0\le2^N$. In this case
\be
{\widehat\phi}(\xi)\chi_2\lll\frac{l\xi}{2^N}\rrr
=m_{l/2^N}\Big(\frac{\xi}{2^{N}}\Big){\widehat\phi}(2\xi),
\ \ \ \ l=0,1\ddd 2^N-1,
\label{33}
\ee
where $m_{l/2^N}(\xi)=\chi_2(l\xi)m_0(\xi)$, $\deg m_{l/2^N}<2^{N+1}$.
It is clear that a wavelet mask $n_0$ can also be chosen in such a way that
$\deg n_0\le2^N$, and we have
\be
{\widehat\psi}(\xi)\chi_2\lll\frac{l\xi}{2^N}\rrr
=n_{l/2^N}\Big(\frac{\xi}{2^{N}}\Big){\widehat\phi}(2\xi),
\ \ \ \ l=0,1\ddd 2^N-1,
\label{34}
\ee
where $n_{l/2^N}(\xi)=\chi_2(l\xi)n_0(\xi)$, $\deg n_{l/2^N}<2^{N+1}$.
Taking the Fourier transform of ~(\ref{33}),~(\ref{34}), we see that
the matrix of the system~(\ref{31}),~(\ref{32}) looks as follows:

$$
\lll
\begin{array}{llllllll}
 g_0 & g_1&\dots& g_{2^N-1}& g_{2^N}&0&\dots&0
 \\
 0& g_0 &\dots& g_{2^{N}-2}& g_{2^{N}-1}& g_{2^N}&\dots&0
 \\
 \hdots&\hdots&\hdots&\hdots&\hdots&\hdots&\hdots&\hdots
 \\
  0& 0 &\dots& g_0& g_1& g_2&\dots&g_{2^N}
  \\
  h_0 & h_1&\dots& h_{2^N-1}& h_{2^N}&0&\dots&0
 \\
 0& h_0 &\dots& h_{2^{N}-2}& h_{2^{N}-1}& h_{2^N}&\dots&0
 \\
 \hdots&\hdots&\hdots&\hdots&\hdots&\hdots&\hdots&\hdots
 \\
  0& 0 &\dots& h_0&  h_1& h_2&\dots&h_{2^N}
\end{array}
\rrr
$$
The determinant of this matrix is so called resultant.
The resultant is not equal to zero if and only if the algebraic
polynomials with the coefficients $g_0,g_1\ddd g_{2^N}$ and
$h_0,h_1\ddd h_{2^N}$ respectively do not have joint zeros
 (see, e.g., \cite{L}).   But this holds
 because the trigonometric polynomials
 $m_0$ and $n_0$ do not have joint zeros by
 construction (taking care of not adding extra zeros).

\bibliographystyle{amsplain}

\end{document}